\documentclass[12pt]{article}
\usepackage{latexsym}
\title{NORMAL SUBGROUPS OF DIFFEOMORPHISM AND HOMEOMORPHISM GROUPS OF $R^n$ AND OTHER OPEN MANIFOLDS}
\author{Paul A. Schweitzer, S.J. \\Pontif\'\i cia Universidade Cat\'olica do Rio de
Janeiro
\\paul37sj@gmail.com}
\usepackage{amssymb}
\usepackage{amsmath}
\usepackage{graphics}

\newtheorem{thm}{Theorem}[section]
\newcommand{\bthm}{\begin{thm}}
\newcommand{\ethm}{\end{thm}}
\newtheorem{prop}[thm]{Proposition}
\newcommand{\bpr}{\begin{prop}}
\newcommand{\epr}{\end{prop}}
\newtheorem{defn}[thm]{Definition}
\newcommand{\bdf}{\begin{defn}}
\newcommand{\edf}{\end{defn}}
\newtheorem{lem}[thm]{Lemma}
\newcommand{\blm}{\begin{lem}}
\newcommand{\elm}{\end{lem}}
\newtheorem{cor}[thm]{Corollary}
\newcommand{\bcr}{\begin{cor}}
\newcommand{\ecr}{\end{cor}}
\newtheorem{conj}[thm]{Conjecture}
\newcommand{\bcj}{\begin{conj}}
\newcommand{\ecj}{\end{conj}}
\newtheorem{rmk}[thm]{Remark}
\newcommand{\brk}{\begin{rmk}}
\newcommand{\erk}{\end{rmk}}
\newtheorem{ex}[thm]{Example}
\newcommand{\bex}{\begin{ex}}
\newcommand{\eex}{\end{ex}}

\newcommand{\sm}{\smallsetminus}
\newcommand{\supp}{{\rm Supp}}
\newcommand{\diff}{{\rm Diff}}

\newcommand{\Int}{{\rm Int}}
\newcommand{\I}{{\rm I}}
\newcommand{\Ker}{{\rm Ker}}
\newcommand{\qed}{\hfill$\Box$}

\begin{document}
\maketitle

\begin{abstract}

We determine all the normal subgroups of the group of $C^r$
diffeomorphisms of ${\mathbb R}^n$, $1\leq r\leq \infty$, except
when $r=n+1$ or $n=4$, and also of the group of homeomorphisms of ${\mathbb R}^n$ ($r=0$). We
also study the group $A_0$ of diffeomorphisms of an open manifold
$M$ that are isotopic to the identity. If $M$ is the interior of a
compact manifold with nonempty boundary, then the quotient of
$A_0$ by the normal subgroup of diffeomorphisms that coincide with
the identity near to a given end $e$ of $M$ is simple. (This version corrects
an error involving the diffeomorphisms of ${\mathbb R}^1$.)

\noindent{\bf keywords:} diffeomorphism groups, homeomorphism
groups, simple groups, normal subgroups, open manifolds

\medskip

\end{abstract}

\section{Introduction}\label{intro}

Let $A_0= A^r_0({\mathbb R}^n)=\diff^r_0({\mathbb R}^n)$ be the
group of $C^r$ diffeomorphisms of ${\mathbb
R}^n$ that are $C^r$ isotopic to the identity, $1\leq r\leq \infty$, or the homeomorphisms isotopic to the identity, if $r=0$. (For simplicity, we shall also refer to
homeomorphisms as diffeomorphisms of class $C^0$.)
It is well known that these groups coincide with the groups of $C^r$ diffeomorphisms of ${\mathbb
R}^n$ that preserve the orientation. Let $A_K$ be
the subgroup of all diffeomorphisms that have compact support, and let $A_{I}$
be the subgroup of those that are isotopic to the identity by a
compactly supported isotopy.
Denote the group of orientation preserving diffeomorphisms of $S^n$ by $\diff^r_+(S^n)$ and let
$\kappa: A_K\to \pi_0\diff^r_+(S^n)$ be
the homomorphism that takes a diffeomorphism $f\in A_K$ to the $C^r$ isotopy class of its extension to $S^n={\mathbb
R}^n\cup\{\infty\}$. In particular, $A_{I} \subseteq \Ker(\kappa)$.
The group $\pi_0\diff^r_+(S^n)$ is easily seen to be abelian (by taking representative diffeomorphisms supported on disjoint balls),
and if $n\neq 4$ and $r\geq 1$ it is known that
it is isomorphic to $\theta_{n+1}$, the finite abelian Kervaire-Milnor group of orientation preserving diffeomorphisms of
$S^n$ up to h-cobordism \cite{KM}. Our main theorem is the following (with a stronger version in Theorem \ref{nonordiff} below).

\bthm The groups $\kappa^{-1}(\eta)$, where $\eta$ varies over all
the subgroups of $\pi_0\diff^r_+(S^n)$, are normal subgroups of
$A_0^r({\mathbb R}^n)$ if $r\geq 1$. If $n>1$ and $r\neq n+1$, they are the
only proper normal subgroups. When $n=1$, there are three proper normal subgroups, $A_{I}=A_K$ and the
diffeomorphisms that coincide with the identity near $+\infty$ or near $-\infty$.
For homeomorphisms ($r=0$),
$A_{I}=A_K$ is the only proper normal subgroup.
\label{mainthm}\ethm

\noindent{\bf Correction.} In the published version of this paper in {\it Ergodic Theory and Dynamical Systems} 31 (2011), pp. 1835-1847, the statements of Theorems \ref{mainthm} and \ref{nonordiff}, Corollary \ref{cornonor}, and the paragraph after Theorem \ref{mainthm} containing an application, were not correct in the case that the dimension is $n=1$, because the real line ${\mathbb R}$ has two ends. The corrected statements are given in this version.

As an application, note that if $f$ is any orientation preserving diffeomorphism of ${\mathbb R}^n$ (of class $C^r$, $r\geq 0$, with $n>1$) that does not have compact support---for example, any affine diffeomorphism that is not the identity---then every orientation preserving diffeomorphism can be expressed as a composition of conjugates in $A_0$ of $f$ and $f^{-1}$. When $n=1$, every orientation preserving diffeomorphism can be expressed as a composition of conjugates in $A$ of $f$ and $f^{-1}$.

In the 1960's and early 1970's, algebraic properties of
homeomorphism and diffeomorphism groups were studied extensively
in work of Anderson \cite{And}, Epstein \cite{Eps}, Herman
\cite{H}, Thurston \cite{Th}, and Mather \cite{Ma1, Ma2}, among
others. The fundamental result is Mather's theorem \cite{Ma1, Ma2,
B} that the group of $C^r$ diffeomorphisms of a smooth manifold
$M$ isotopic to the identity by isotopies with compact support is
a simple group, provided that $r\neq \dim(M)+1$. This theorem had
been proven by Thurston \cite{Th} when $r=\infty$. (The case when
$r=0$ is fairly easy, and the case when $r=\dim(M)+1$ remains
open.) These results are based on Epstein's theorem \cite{Eps}
that under certain very general hypotheses, the commutator
subgroup of a group of homeomorphisms is simple. A survey of
diffeomorphism groups can be found in Banyaga \cite{B} and some
recent results appear in the Morita symposium \cite{P}.

In my lecture at the AMS Symposium on Pure Mathematics at Stanford
University in August 1976, I presented theorems about groups of
diffeomorphisms of ${\mathbb R}^n$. After the lecture, Wensor
Ling, a graduate student working on his doctorate under Thurston
at MIT and Princeton, mentioned that he had obtained similar
results in his doctoral research. I decided to leave the
publication of the results to the future Dr. Ling. Two years
later, Dusa McDuff published a study \cite{McD} of the lattice of normal subgroups of the group of diffeomorphisms of an open
manifold $P$ that is the interior of a compact manifold with
nonempty boundary. She referred to Ling's preprint \cite{L2} and
unpublished work of mine for the proof of a fundamental lemma
stating that $A_0(P)/G_i$ (in her notation, see \cite{McD} p. 353)
is a simple group. Since the proof of this result and the paper \cite{L2}
have apparently never been published, I include its proof here as
Theorem \ref{mainlemma}. I also prove some new results for any
open manifold $M$ of dimension $n$. All the results are valid for diffeomorphism
groups of class $C^r$ (or $C^{r,\alpha}$ in the sense of Mather \cite{Ma1, Ma2}), $r=1,2,\dots,\infty$ and also for
homeomorphism groups ($r=0$), excluding $C^{n+1}$ diffeomorphisms (and with restrictions on $\alpha$
for $C^{r,\alpha}$ diffeomorphisms when $r=n$ or $n+1$) in
Theorems \ref{mainthm} and \ref{nonordiff}, since their proofs use Mather's theorem. The results involving ${\mathbb R}^n$ are proven in Section \ref{Rnproofs} using the results on general open manifolds, which are proven in Section \ref{proofs}.
I would like
to thank Tomasz Rybicki for encouraging me to publish this paper
after all these years.

\section{Definitions and Results}\label{results}

We recall some definitions from McDuff's paper \cite{McD}. Let $M$
be a connected $n$-manifold. For $r=1,2,\dots,\infty$, we set
$A=A(M)= \diff^r(M)$, the group of all $C^r$ diffeomorphisms of
$M$ with the $C^r$ compact-open topology, and for $r=0$ we let
$A={\rm Homeo(M)}$ be the group of homeomorphisms of $M$. Let
$A_0$ be the identity component of $A$, so that it consists of the
diffeomorphisms that are $C^r$ isotopic to the identity map $\I$
(or the homeomorphisms that are topologically isotopic to $\I$,
if $r=0$). The support $\supp(g)$ of a diffeomorphism $g\in
A$ is the closure of the set $\{x\in M\ |\ g(x)\neq x\}$, and the
support of an isotopy $\{g_t\ |\ t\in [0,1]\}$ with $g_0=\I$ is
the closure of the union of the supports $\cup_{t\in [0,1]}\supp(g_t)$.
All manifolds are connected.

We may also consider diffeomorphisms of class $C^{r,\alpha}$ where $\alpha:[0,\epsilon]\to {\mathbb R}$ is a  {\em modulus of continuity} (see \cite{Ma1} for the definitions). (In particular, the $r$th derivatives will be H\"older continuous with exponent $k\in (0,1)$ if $\alpha(t)=t^k$, or Lipschitz if $\alpha(t)=t$.)  A mapping from an open set in ${\mathbb R}^m$ to ${\mathbb R}^n$ is $C^{r,\alpha}$ if its $r$th derivative exists and is {\em locally $\alpha$-continuous}. For simplicity of notation, we also denote
$C^r$ diffeomorphisms by $C^{r,\alpha}$, understanding in this case that $\alpha=\emptyset$.
A $C^{r+1}$ mapping is $C^{r,\alpha}$ for every modulus of continuity $\alpha$.
Mather observes that sums, products, compositions and inverses of $C^{r,\alpha}$ mappings with $r\geq 1$ are $C^{r,\alpha}$. That is enough for all of our constructions to work, so we make the following convention.

\vskip 0.3cm
\noindent {\bf Convention.}
Diffeomorphisms will be of class $C^{r,\alpha}$ with $r\geq 1$ (including $C^r=C^{r,\emptyset}$) or else $C^0$ (i.e., homeomorphisms). Isotopies are required to have the same differentiability as the diffeomorphisms.

\bthm Let $h\in A_0$ be a diffeomorphism isotopic to the identity
map of an open (i.e., non-compact) manifold $M$ such that $\supp(h)$ is not compact,
and let $H$ be any normal subgroup of $A_0$ containing $h$. Then
the subgroup $A_I$ of $A_0$ consisting of all diffeomorphisms that are
isotopic to the identity map of $M$ by isotopies with compact
support is contained in $H$. \label{thm2} \ethm

In the case that the manifold $M$ is ${\mathbb R}^n$, we have the following
stronger version of Theorem \ref{mainthm}, using the same notation. Let $A=\diff^{r,\alpha}({\mathbb R}^n)$ be the group of all $C^{r,\alpha}$ diffeomorphisms of ${\mathbb R}^n$, $1\leq r<\infty$ with $\alpha$ a modulus of continuity, or else
$A=\diff^r({\mathbb R}^n)$ with $0\leq r\leq\infty$. Let $A_0$
be the subgroup of diffeomorphisms isotopic to the identity, and $A_K$ those with compact support.
In view of the statements of Mather's theorems \cite{Ma1, Ma2}, we say that a pair
$(r,\alpha)$ where $1\leq r<\infty$ and $\alpha$ is a modulus of continuity, {\em satisfies Mather's condition} if one of the following holds:
\begin{enumerate}
\item $1\leq r<n$ or $r>n+1$, and $\alpha$ is arbitrary;
\item $r=n$, $\alpha:[0,\infty)\to
[0,\infty)$, and there is a positive number $\beta<1$ such that
$\alpha(tx)\leq t^\beta\alpha(x)$ for all $x\geq 0$ and all $t\geq 1$;
\item $r=n+1$ and $\alpha(x)=x^\beta$ for some positive number $\beta\leq 1$.

\end{enumerate}

\bthm The subgroup $A_0$ of $A({\mathbb R}^n)$ consisting of
diffeomorphisms of ${\mathbb R}^n$ isotopic to $I$ (which
coincides with the group of orientation preserving
diffeomorphisms),
 and the subgroups $\kappa^{-1}(\eta)$ where $r\geq 1$ and $\eta$ varies over all the subgroups
 of $\pi_0\diff^r_+(S^n)$ (or, equivalently, of $\theta_{n+1}$, if $n\neq 4$)
 are normal subgroups of $A_0$ and (except possibly if $n=4$) also of $A$.
 There are no other
 proper normal subgroups of $A_K$, of $A_0$, and of $A$, for $C^r$ diffeomorphisms with $1\leq r\neq n+1$, and also for $C^{r,\alpha}$ diffeomorphisms if $(r,\alpha)$ satisfies Mather's condition,
 provided that $n>1$. When $n=1$, there are no other proper normal
 subgroups of $A_K$ and $A$, but $A_0$ has two additional
 proper normal subgroups, consisting of the diffeomorphisms that
 coincide with the identity near to $+\infty$ or near to $-\infty$.
 For $r=0$, the only proper normal subgroups of $A_0$ and $A$ are
 $A_K=A_I$ and $A_0$ in $A$. \label{nonordiff}
\ethm

\bcr Suppose $n>1$. When $(r,\alpha)$ satisfies Mather's condition or when $1\leq r\neq n+1$ and $\alpha=\emptyset$, the only subgroups that appear in any descending normal series of subgroups starting with $A=\diff^{r,\alpha}({\mathbb R}^n)$, $A_0=\diff^{r,\alpha}_0({\mathbb R}^n)$, or $A_K^{r,\alpha}({\mathbb R}^n)$, besides these groups themselves and the trivial subgroup, are the subgroups $\kappa^{-1}(\eta)$, where $\eta$ varies over all subgroups of $\pi_0\diff^1_+(S^n)$. \label{cornonor}
\ecr

The following result, used by McDuff \cite{McD} to study the
lattice of normal subgroups of the group of diffeomorphisms
isotopic to the identity on the interior $P$ of a compact manifold
$\bar P$ with non-empty boundary, is fundamental, and will be
proven in Section \ref{proofs}. Say that an end $e$ of an open
manifold $M$ is a {\em product end} if it has a basis of
neighborhoods homeomorphic to $N\times (0,1/k]$ for some fixed
closed $(n-1)$-manifold $N$. Such neighborhoods are called {\em product neighborhoods} of $e$. Let $A_e$ be the subgroup of $A_0$
consisting of diffeomorphisms that are the identity in a
neighborhood of $e$.

\bthm If $e$ is a product end of $M$ then $A_0/A_e$ is a simple group.
 \label{mainlemma}\ethm

In the notation of \cite{McD}, this says that if $e$ is an end of
$P=\Int(\bar P)$ and $G_i=A_e$, then the group $A_0/G_i$ is
simple.

The four following propositions and their corollaries are steps in
the proofs of the theorems. Every ball considered in an
$n$-manifold $M$ will be a {\em closed collared $n$-ball}, i.e.,
the image of the closed unit ball in $\mathbb{R}^n$ under a $C^r$
embedding of $\mathbb{R}^n$ in $M$. Let $B_m, m=1,2,\dots,$ be a
sequence of disjoint balls in an open manifold $M$ that converge
to an end $e$ of $M$, i.e., every neighborhood of $e$ contains all
but a finite number of the balls $B_m$.

\bpr Every diffeomorphism of $M$ supported in $\cup_{m=1}^\infty B_m$ is isotopic to the identity map $\I$.
\label{infballs} \epr

In particular, every diffeomorphism of an open manifold supported
on a finite set of disjoint balls belongs to $A_0$.

\bpr Let $h\in A_0$ be a diffeomorphism of $M$ whose support meets
every neighborhood of the end $e$ and let $H$ be a normal subgroup
of $A_0$ containing $h$. If $\dim(M)\geq 2$, then there exists
$f\in H$ such that $f(B_m)= B_{m+1}$ for every $m$. Furthermore,
every diffeomorphism of $M$ supported in $\cup_{m=1}^\infty B_m$
belongs to $H$.
If $\dim(M)=1$ (so that
$M\approx {\mathbb R}$), the same conclusions hold, provided that the disjoint intervals
$B_m$ are in a monotone increasing (or decreasing) order on the line. \label{disjointballs} \epr

\bpr Let $e$ be a product end of an open manifold $M$, and let $H$
be a normal subgroup of $A_0=A_0(M)$ containing a diffeomorphism
$h\in A_0$ whose support $\supp(h)$ meets every neighborhood of
$e$. Then there exist $f_0\in H$ and a closed product neighborhood
$N_0$ of $e$ such that the sets $f_0^m(N_0)$ for $m\geq 0$ are
nested ($f_0^{m+1}(N_0)\subset \Int \ f_0^m(N_0)$) and form a
neighborhood base for $e$. \label{nbdbase} \epr

\noindent Note that the diffeomorphism $f_0$ is related to the
translations of \cite{L1}.

\bcr Under the same hypotheses on $h$ and $H$, every diffeomorphism $g$ of $M$ with support
contained in a product neighborhood of $e$
belongs to $H$. \label{prodsupp}\ecr

The following Factorization Lemma is an essential step in the
proofs. It is well-known, but since there are a variety of
versions, we give explicitly the version we shall use.

\bpr {\bf Factorization Lemma.} Given open sets $U_1,\dots,U_k$ in
a manifold $M$, a compact set $C$, and a $C^r$ isotopy $F:M\times
[0,1]\rightarrow M$ of the identity map $\I$ to a diffeomorphism
$f$ of $M$ such that $F(C\times [0,1]) \subset  \cup_{j=1}^k U_j$,
there are finitely many diffeomorphisms $g_0=\I,g_1,\dots, g_p\in
A_0$, a function $j:\{1,\dots,p\}\rightarrow \{1,\dots,k\}$,
isotopies $F_i$ of $\I$ to $f_i= g_ig_{i-1}^{-1}$ with
$\supp(F_i)\subseteq U_{j(i)}$, and an isotopy of $\I$ to
$fg_p^{-1}$ supported on $M\sm C$. In particular, $g_p|_C=f|_C$.
\label{factorization}\epr

In other words, the given isotopy can be factored into `small'
isotopies on $C$, each supported in one of the given open sets,
plus one other isotopy supported in the complement of $C$. The
idea of the proof is as follows. The proof of Lemma 3.1 of
\cite{PS} for $C^r$ diffeomorphisms ($r\geq 1$) shows that each
$t\in [0,1]$ has a neighborhood $[t_1,t_2]\subseteq [0,1]$ such
that the isotopy $F\circ(F_{t_1}^{-1}\times \I)$ restricted to
$M\times [t_1,t_2]$ has such a factorization on $F_{t_1}(C)$, and
Theorem 5.1 of \cite{EK} implies the same result when $r=0$. Then
compactness of the interval gives the desired conclusion. The
following corollary is an easy consequence.

\bcr Suppose that an isotopy from the identity map of $M$ to a diffeomorphism $f\in A_0$
has compact support contained in an open set $U$. Then the given isotopy can be
factored into a finite composition of isotopies, each with support in an open ball in $U$.
\label{factorization2} \ecr

It is clear from Theorem \ref{mainthm} that if $r\geq 1$ there exist $C^r$
diffeomorphisms of ${\mathbb R}^n$ in many dimensions that have compact support and are isotopic to the
identity $I$, but not by isotopies with compact support. Any isotopy
of such a diffeomorphism $f$ to $I$ cannot be factored into a finite
composition of isotopies supported in open balls. On the other
hand, when $r=0$ it is well-known that every homeomorphism of
${\mathbb R}^n$ with compact support is isotopic to the identity
map by a compactly supported topological isotopy. This is easy to
see by the `Alexander trick' of shrinking the support into the
origin.

Other examples of diffeomorphisms isotopic to $I$, but not by
isotopies with compact support, are the following.

\bex Diffeomorphisms of $M=N\times {\mathbb R}$ in $A_0(M)$
isotopic to the identity through diffeomorphisms with compact
support, but not by compactly supported isotopies.  \eex

Let $M=N\times {\mathbb R}$ where $N$ is a closed manifold and
consider a map $F:N\times {\mathbb R}\rightarrow N$ of class $C^r,
0\leq r\leq\infty$, such that the restriction of $F$ to $N\times
[0,1]$ is an isotopy of the identity to itself and such that
$F_t=\I$ when $t\in {\mathbb R}\sm (0,1)$. Suppose that for some
point $x_0\in N$ (and hence for every $x\in N$) the loop $\pi
\circ \gamma$ represents a non-zero element of the fundamental
group $\pi_1(N)$, where $\pi: N\times [0,1]\rightarrow N$ is the
projection and $\gamma(t)=F(x_0,t)$, $t\in [0,1]$. Then the diffeomorphism $f$
of $N\times {\mathbb R}$ with compact support defined
$f(x,t)=(F(x,t),t)$ is isotopic to $\I$, since the support of $f$
can be pushed off to infinity, but no isotopy of $f$ to $\I$ can
have compact support. Furthermore, $f$ cannot be the composition
of a finite number of diffeomorphisms supported in closed balls.

Concretely, in dimension $n=2$, let $M=S^1\times {\mathbb R}$ and
let $f(z,t)=(e^{2\pi i \lambda(t)z},t)$ where $\lambda: {\mathbb
R}\rightarrow [0,1]$ is a smooth function that is $0$ for $t\leq
0$ and $m\neq 0$ for $t\geq 1$ (an $m$-fold Dehn twist). In
dimension $n>2$ this example can be multiplied by a closed
$(n-2)$-manifold $N'$ to give such a diffeomorphism of $M=N'\times
S^1\times {\mathbb R}$.

\section{Proofs of results on diffeomorphisms of ${\mathbb R}^n$}\label{Rnproofs}

In this section we prove Theorems \ref{mainthm} and \ref{nonordiff}
 about diffeomorphisms of
${\mathbb R}^n$. We suppose the other results that are mentioned
above and will be proven in the last section. The following
theorem of Mather will also be used. Recall Mather's
condition on the pair $(r,\alpha)$ as defined before the statement of Theorem \ref{nonordiff}.

\bthm (Mather \cite{Ma1, Ma2, B}) For any connected $n$-dimensional smooth manifold
$M$, the group $A_I(M)$ of $C^r$ diffeomorphisms of $M$ that are
isotopic to the identity $I$ by isotopies with compact support is simple, if $1\leq r\leq \infty$ with $r\neq n+1$.
The same holds for $C^{r,\alpha}$ diffeomorphisms if $(r,\alpha)$ satisfies Mather's condition.
 \label{Mather}
\ethm

Now, as in the Introduction, we consider the homomorphism $\kappa: A_K\to \pi_0\diff_+(S^n)$ that takes $g\in A_K$ to the isotopy class of its extension $\bar g$ over $S^n$. The group $\pi_0\diff_+(S^n)$, as observed in the introduction,
is abelian, independent of the class of
differentiability $C^r$ provided that $r\geq 1$, and also, if
$n\neq 4$, isomorphic to $\theta_{n+1}$ and finite.

\bpr If $g\in A_K=A_K({\mathbb R}^n)$ satisfies $\kappa(g)=0\in
\pi_0\diff_+(S^n)$, then $g\in A_I$. Consequently, if two
diffeomorphisms (of differentiability at least $C^1$) of
${\mathbb R}^n$ with compact support have
smoothly isotopic extensions on the sphere $S^n={\mathbb
R}^n\cup\{\infty\}$ (of class $C^1$), then they are isotopic on
${\mathbb R}^n$ by a smooth isotopy with compact support.
\label{undoisot}\epr It is enough to suppose that a smooth isotopy
is $C^1$ since a $C^1$ isotopy between $C^r$ (or $C^{r,\alpha}$)
diffeomorphisms can easily be smoothed to an isotopy of the same
differentiability class. The Proposition also holds for homeomorphisms,
but we shall not need that fact. The proof of the Proposition
uses the following lemma.

\blm Let $U\subseteq {\mathbb R}^n$ be an open set that contains
the origin $0$ and let $s\to F_{t,s}: U\to F_{t,s}(U),\ s\in
[0,1]$, be a smooth family of isotopies into ${\mathbb R}^n$
parameterized by $t\in N$, where $N$ is a compact smooth manifold with
boundary, such that $F_{t,s}(0)=0$ for every $t$ and $s$. Then
there exist $\epsilon>0$ and a smooth family $H_{t,s}$ of
isotopies of ${\mathbb R}^n$ with support in a closed ball
$B_d(0)\subseteq U$ centered at $0$ such that, for
every $t$, $H_{t,0}=I$, $H_{t,1}$ coincides with $F_{t,1}\circ
F_{t,0}^{-1}$ on the ball $B_\epsilon(0)$, and $H_{t,s}=I$
whenever $F_{t,s}=I$ on $U$ for all $s$. \label{blowup}\elm

\noindent{\bf Proof.} Note that this is an extension of the
Factorization Lemma \ref{factorization} in two ways: first, it
involves a smoothly parameterized family, and second, it does not require that the isotopies be surjective.

Suppose we are given a smooth family $F:U\times N\times [0,1]\to
{\mathbb R}^n$ of isotopies into ${\mathbb R}^n$ such that
$F(0,t,s) = 0$ for every $t$ and $s$. Set $F_{t,s}(x) =
F(x,t,s)$. Note that $t\in N$ is the parameter and $s\in [0,1]$ is
the index of isotopy. Since $0$ is preserved, if $S_\epsilon$ is a small $(n-1)$-sphere centered at $0$ and contained in $U$, the set $\cup_{t,s}F_{t,s}(S_\epsilon)$ is compact and does not contain $0$, so
there is a closed ball $B_d=B_d(0)$ centered at $0$ with positive radius $d$ contained in the intersection of the sets $F_{t,s}(U)$ for all $(t,s)\in N\times[0,1]$. Let $B_{d'}$ be a smaller concentric ball in the interior of $B_d$ and let  $\lambda:{\mathbb R}^n\to [0,1]$ be a smooth
function such that $\lambda=0$ on $B_{d'}$ and
$\lambda=1$ on a neighborhood of ${\mathbb R}^n\sm \Int(B_d)$.
Set $G_{t,s,s'}=F_{t,s}\circ F_{t,s'}^{-1}$ and note that $G_{t,s,s'}$ is defined on all of $B_d$. Hence, given $s,s'\in [0,1]$ and $t\in N$, we may define a mapping $\bar G_{t,s,s'}: B_d \to {\mathbb R}^n$ depending on $t, s$ and $s'$ by setting
$$\bar G_{t,s,s'}(x)= \lambda(x)x + (1-\lambda(x))G_{t,s,s'}(x).$$
The derivative of $\bar G_{t,s,s'}$ with
respect to the variable $x$ is
\begin{equation}
\begin{array}{rl}
\!D\bar G_{t,s,s'}(x)(v)=&\!\!\nabla\lambda(x)\!\cdot\! v\ x + \lambda(x)v
\\&-\nabla\lambda(x)\!\cdot\! v\ G_{t,s,s'}(x) +(1-\lambda(x))DG_{t,s,s'}(x)(v) \nonumber\\
=&\!v + \nabla\lambda(x)\!\cdot\! v(x-G_{t,s,s'}(x))
\\&+(1-\lambda(x))(DG_{t,s,s'}(x)(v)-v),\nonumber\\
\end{array}
\end{equation}
where $v\in {\mathbb R}^n$. Now the gradient $\nabla\lambda(x)$ and the function
$1-\lambda(x)$ are bounded, $F$ is smooth of at least class $C^1$,
and $\bar G_{t,s,s'}=I$ when $s=s'$, so by uniform continuity of $\bar G_{t,s,s'}$ and its
derivative $D\bar G_{t,s,s'}$ as functions of $x, t, s,$ and $s'$, there
exists $\delta>0$ such that whenever $|s-s'|\leq\delta$ the
second and third terms in the formula for the derivative will be sufficiently small so that $D\bar G_{t,s,s'}$
will be invertible for all $t\in N$. Hence $\bar G_{t,s,s'}$ is a local diffeomorphism on
$B_d$, and since it coincides with the identity $I$ outside a
compact set in $\Int(B_d)$, it is a diffeomorphism of $B_d$ and extends
to a global diffeomorphism of ${\mathbb R}^n$ with support inside
$B_d$. Thus if $s'<s''\leq s'+\delta$, the mapping
$s\in[s',s'']\mapsto \bar G_{t,s,s'}$ defines a smooth family of
isotopies of ${\mathbb R}^n$ with support in $B_d$, such that
$\bar G_{t,s',s'}=I$ and on $B_{d'}$ the restrictions of
$\bar G_{t,s'',s'}$ and $F_{t,s''}\circ F_{t,s'}^{-1}$ agree.

Now take a partition $0=s_0<s_1<\dots<s_k=1$ of the interval
$[0,1]$ such that each subinterval $[s_{i-1},s_i]$ has length at
most $\delta$. The preceding construction with $s\in [s',s'']=
[s_{i-1},s_i]$ provides a smooth family of isotopies supported on
$B_d$ transforming $I$ into a diffeomorphism $g_{t,i},\
i=1,\dots,k, $ that agrees with $F_{t,s_i}\circ F_{t,s_{i-1}}^{-1}$
on $B_{d'}$. Then the composition $g_{t,k}\circ\dots
g_{t,2}\circ g_{t,1}$
 agrees with $F_{t,1}\circ F_{t,0}^{-1}$ on a possibly smaller ball
 $B_\epsilon(0)$, so the composition of these
 families of isotopies gives the required smooth family of isotopies
 $H_{t,s}$. It is easy to check from the definition of $\bar G$ that
 $H_{t,s}=I$ whenever $F_{t,s}=I$ on $U$ for all $s$.\qed

\vskip 0.5cm

\noindent{\bf Proof of Proposition \ref{undoisot}.} If $r=0$, the
Alexander trick shows that $A_K=A_I$, so suppose that $r\geq 1$.
Let $g\in A_K$ satisfy $\kappa(g)=0\in \pi_0\diff_+(S^n)$. Then
there is a $C^1$ isotopy on $S^n= {\mathbb R}^n\cup \{\infty\}$
from $I$ to the extension $\bar g$ over $S^n$. We shall modify
this isotopy $F:S^n\times [0,1]\to S^n$ so that it will have
compact support contained in ${\mathbb R}^n$, and then it will
restrict to the desired isotopy of $I$ to $g$.

Identify $S^n$ with ${\mathbb R}^n\cup\{\infty\}$ by stereographic
projection through the point $\infty$, which is taken to be the
unit positive vector in ${\mathbb R}^{n+1}$ on the axis orthogonal
to ${\mathbb R}^n$. Look at the loop $t\mapsto F(\infty,t)$.
Modify $F$ so that this loop becomes the constant loop at
$\infty$. To do this when $n\geq 2$, first perturb $F$ slightly so
that the curve does not pass though $0\in {\mathbb R}^n$, and
then modify each $F_t:=F(\cdot,t)$ by the rotation in $SO(n+1)$
that takes $F_t(\infty)$ to $\infty$ while fixing the
$(n-1)$-plane orthogonal to both $F_t(\infty)$ and $\infty$ when
they are distinct. In case $n=1$, it may be necessary to follow
the isotopy $t\mapsto F_t$ by a rotation of $S^1$ so that the loop
becomes null-homotopic.

Next, the Gram-Schmidt orthogonalization process applied to the
derivatives $D_\infty F_t$ on the tangent space $T_\infty S^n$
gives canonical smooth paths in $GL(n,{\mathbb R})$ transforming
each $D_\infty F_t$ into an element of $O(n)$. Projecting these
paths in $GL(n,{\mathbb R})$ down onto $S^n$ and composing them
with the diffeomorphisms $F_t$ defines a smooth family of
isotopies on a ball in $S^n$ containing $\infty$, and we identify its interior with an open set $U\subset{\mathbb R}^n$ so that $\infty$
is identified with $0\in U$. Applying Lemma \ref{blowup} with
$N=[0,1]$ and then projecting the resulting diffeomorphisms to
$S^n$, we obtain a family of isotopies $H_{t,s}$ supported in an
open set containing $\infty$. Replace $F_t$ by $H_{t,1}\circ F_t$,
so that the derivatives $D_\infty F_t$ of the new family $F_t$
will be orthogonal for each $t$. The loop $t\mapsto D_\infty F_t$
in $SO(n)$ can be made null-homotopic, following it by its inverse
loop, if necessary. Then we may shrink it by composing with
elements of $SO(n)$ acting on $S^n$ so that $D_\infty F_t$ becomes
the identity map $I$ for all $t$. We approximate $F:S^n\times [0,1]\to S^n$ so that it becomes $C^\infty$ in a neighborhood of $\{\infty\}\times [0,1]$. All these steps can be carried
out so that $F_0=I$ and $F_1=\bar g$ are not changed.

Now, still denoting the family of diffeomorphisms by $F_t$,  there is a smooth family of isotopies from $I$ to $F_t$ on a
coordinate chart around $\infty\in S^n$ with $\infty$ identified
with $0\in {\mathbb R}^n$, defined as follows:

\begin{equation} \tilde F(x,t,s) =
\begin{cases}
(x,t)& \text{if}\ s=0,\\
    (s^{-1}F_t(sx),t) & \text{if}\ s\in (0,1],
    \end{cases}
    \end{equation}
    so that $\tilde F(\cdot,t,1)=F_t$.
It is easy to check that $\tilde F$ has the same class of differentiability as $F$, by using the Taylor expansions of the $F_t$'s, since we have made $F$ to be $C^\infty$ near $\{0\}\times [0,1]$.

Applying Lemma \ref{blowup} again to this family of isotopies in a
small coordinate chart that takes $\infty$ to $0\in {\mathbb
R}^n$, we may replace $F_t$ by $H_t^{-1}\circ F_t$ in the
coordinate chart. The result is a new smooth family of
diffeomorphisms $F_t, t\in [0,1]$, that coincide with the identity
in a fixed ball around $\infty$ and still satisfy $F_0=I$ and
$F_1=\bar g$. Thus they restrict to a smooth isotopy with compact
support from $I$ to $g$ on ${\mathbb R}^n$, as desired.

Now if two diffeomorphisms $f$ and $g$ of ${\mathbb R}^n$ with compact support have
isotopic extensions $\bar f$ and $\bar g$ on the sphere $S^n$,
then there is an isotopy of $I$ to $\bar f^{-1}\circ\bar g$, so by
the first part of the Proposition, $f^{-1}\circ g$ is isotopic to
$I$ by an isotopy of ${\mathbb R}^n$ with compact support, which
implies the second conclusion. \qed

\bpr Let $g\in A_K=A_K({\mathbb R}^n)$ and $f\in
A=\diff^r({\mathbb R}^n)$. If $f$ preserves the orientation, then
there is an isotopy with compact support from $fgf^{-1}$to $g$. If
$f$ reverses the orientation and $n\neq 4$, then there is an
isotopy with compact support from $fgf^{-1}$to $g^{-1}$.
\label{conjug} \epr

\noindent{\bf Proof.} First, suppose that $f$ preserves the
orientation, so that $f\in A_0$, and there is an isotopy $F:
{\mathbb R}^n\times [0,1]\to {\mathbb R}^n$ such that $f_0 =I$ and
$f_1=f$, where $f_t=F(\cdot,t)$. If $g\in A_K$, then the set
$K=\supp(g)= \supp(g^{-1})$ is compact, so $K'=F(K\times
[0,1])$ is also compact. It is easy to check that the isotopy
$t\mapsto f_tgf_t^{-1}$ from $g$ to $fgf^{-1}$ has support in
$K'$.

Let $f$ be reflection in the last coordinate, $$f(x_1,\dots,
x_{n-1},x_n) = (x_1,\dots, x_{n-1},-x_n),$$ and suppose that
$\supp(g)\subset \{x_n>0\}\subset {\mathbb R}^n$. If $n\leq 3$ it
is known that $\pi_0\diff_+(S^n)$ is trivial and we are assuming
that $n\neq 4$, so we may suppose that $n\geq 5$. Define $h\in
A_K$ to coincide with $g$ where $x_n\geq 0$ and with $fgf^{-1}$
where $x_n\leq 0$. Now $f$ extends to $\bar f$ on ${\mathbb
R}^{n+1}$, just reversing the sign of the coordinate $x_n$, and
$\bar f$ commutes both with $h$ and with the stereographic
projection through $\infty=(0,\dots,0,1)$, and therefore also with
the extension $\bar h$ of $h$ to $S^n$. Consequently $\bar h$
extends to a diffeomorphism $\tilde h$ on the ball $D^{n+1}$; it
suffices to let $\bar h(x)=y$ and set $$\tilde
h(x_0,x_1,\dots,x_{n-1},tx_n)=(y_0,y_1,\dots,x_{n-1},ty_n)$$ for
$t\in [-1,1]$. This is a $C^r$ diffeomorphism since $\bar h$ is
the identity in a neighborhood of $S^n\cap\{x_n=0\}$. It follows
that $\bar h$ is pseudo-isotopic to $I$.
 Since we are supposing $n\geq 5$, pseudo-isotopy on a simply connected manifold implies isotopy \cite{Cerf}, so $\bar h$ is isotopic to $I$ on $S^n$. Since $\bar h$ is isotopic to the composition $(fgf^{-1})g$, we conclude that $fgf^{-1}$ is isotopic to $g^{-1}$, as desired.

 Now any $g_1\in A_K$ is isotopic to a diffeomorphism $g$ with compact support in $\{x_n>0\}$ and any orientation reversing $f_1\in A$
 can be factored $ff_0$ where $f$ is the reflection used above and $f_0\in A_0$, so if $n\geq 5$ it follows that $f_1g_1f_1^{-1}$ is
 isotopic to $g_1^{-1}$. \qed

\vskip .5cm

It seems likely that $fgf^{-1}$ is also isotopic to $g^{-1}$ when
$n=4$, but the above proof would need the equivalence of
pseudo-isotopy to isotopy, and that equivalence is apparently
unknown in dimension $4$.

\vskip .5cm

\noindent{\bf Proof of Theorems \ref{mainthm} and
\ref{nonordiff}.} Suppose $r\geq 1$. Note that for every subgroup
$\eta\subseteq\pi_0\diff_+(S^n)$, $\kappa^{-1}(\eta)$ is a subgroup
of $A_0$ and of $A=\diff({\mathbb R}^n)$; by Proposition
\ref{conjug}, it is a normal subgroup of $A_0$ and also, if $n\neq
4$, of $A$. In particular, $\kappa^{-1}(0)= A_I$ and
$\kappa^{-1}(\pi_0\diff_+(S^n))=A_K$ are normal. Furthermore, $A_0$
has index $2$ in $A$ and so it is also normal. Thus all of the
subgroups mentioned in Theorems \ref{mainthm} and \ref{nonordiff}
are normal.

When $n=1$ the results are easy to show, so we suppose that $n>1$, and then ${\mathbb R}^n$ has only one end.
We must show that these are the only non-trivial normal subgroups
when either $r\neq n+1$ or $(r,\alpha)$ satisfies Mather's condition. Let $H$ be a non-trivial normal subgroup of
$A_0$. If $H$ contains a diffeomorphism $h$ with non-compact
support, then by Proposition \ref{disjointballs}, $H$ contains $A_K$,
since every compact set in ${\mathbb R}^n$ is contained in a ball.
Now by Theorem \ref{mainlemma}, $A_0/A_K$ is simple (since
$A_K=A_e$), so $H/A_K=A_0/A_K$ and $H=A_0$. If $H$ is a normal
subgroup of $A$ not contained in $A_0$, then we claim that $H\cap
A_0$ contains an element with non-compact support. To see this,
take $h\in H$ with non-compact support, and conjugate $h$ by an
element $f\in A_0$ so that there are points $x_n, n\in {\mathbb
N},$ tending to infinity such that $f^{-1}hf(x_n)\neq h(x_n)$ for
all $n$. Then $h^{-1}f^{-1}hf\in H\cap A_0$ has non-compact
support. It follows that $H\cap A_0$ contains both $A_K$ (as before)
and also $A_0$ (since $A_0/A_K$ is simple). Since $H$ is not
contained in $A_0$, we have $H=A$.

It remains to consider a normal subgroup $H\neq \{I\}$ of
$A_K$. It is easy to see that $H\cap A_I\neq \{I\}$. For example,
conjugate a non-trivial element $g\in H$ by a small perturbation
$f\in A_K$ so that $I\neq fgf^{-1}g^{-1}\in H\cap A_I$. Then by
Mather's Theorem \ref{Mather}, $A_I\cap H=A_I$ and therefore
$A_I\subseteq H$, provided that either $r\neq n+1$ or $(r,\alpha)$ satisfies Mather's condition.

To complete the proof when $r\geq 1$, it suffices to set $\eta=\kappa(H)$ and
show that $\kappa^{-1}(\eta)=H$. Clearly $H\subseteq \kappa^{-1}(\eta)$. Take any $f\in \kappa^{-1}(\eta)$.
Then there must be some $f_1\in H$ such that $\kappa(f)=\kappa(f_1)$.
Thus $g=ff_1^{-1}$ satisfies $\kappa(g)=0\in \pi_0\diff_+(S^n)$.
By Proposition \ref{undoisot}, $g\in A_I\subset H$, so $f=gf_1\in H$. Thus $H\subseteq \kappa^{-1}(\eta)$, so $H=\kappa^{-1}(\eta)$.

Now let $r=0$. The normal subgroups $A_K$ and $A_I$ coincide,
since a homeomorphism of ${\mathbb R}^n$ with support in a ball is
isotopic to the identity by an isotopy with support in the same
ball, by just shrinking the support of the homeomorphism into the
origin (the `Alexander trick'). As when $r\geq 1$, any
normal subgroup of $A_0$ not contained in $A_K$ coincides with $A_0$, and any normal subgroup $H$ of $A$ not contained in $A_0$ coincides with $A$. Hence it suffices to show that any
non-trivial normal subgroup $H$ of $A_K$ must
coincide with $A_K$. Let $h\in H\sm \{I\}$, let $x_0$ be a point
on the frontier of $\supp(h)$, and choose a sequence of disjoint
balls $B_1, B_2, \dots$ in $M={\mathbb R}^n\sm \{x_0\}$ converging
to the end $e$ where $x_0$ was removed (with care to have them in monotone order on the line, if the dimension $n$ is $1$). Consider any $g\in A_K$
with $\supp(g)\subseteq B_1$. Note that the restriction of $h$ to $M$ belongs to
$H':=\{f\in H| f(x_0)=x_0\}$ and moves points arbitrarily close
to $e$. By applying Proposition \ref{disjointballs} to the normal subgroup $H'$ of $A_0(M)$ and the restrictions of $h$ and $g$ to $M$, we find that $g\in H'\subseteq H$.
Now any homeomorphism $g'\in A_K$ can be conjugated by an element
of $A_K$ to some $g\in A_K$ supported in $B_1$, so $g'$ and hence
$g$ belong to $H$. Thus $H$ must be all of $A_K$. \qed

\section{Proofs of the results for open manifolds}\label{proofs}

This section is independent of the previous section. Diffeomorphisms may have any differentiability class $C^r$ with $0\leq r\leq \infty$ or $C^{r,\alpha}$ with $r\geq 1$ and any modulus of continuity $\alpha$.
We shall use the following easy to prove lemma.

\blm 1. If $B$ and $B'$ are balls in a manifold $M$,
$\gamma$ is a path from a point in the interior of $B$ to another
point in the interior of $B'$, and $U$ is an open set containing
the balls and the image of $\gamma$, then there is a diffeomorphism $g\in
A_0$ such that $g(B)= B'$ and $g$ is isotopic to the identity map
by an isotopy supported in $U$.

2. If $B_m, m=1,2,\dots,$ is a locally finite sequence of disjoint balls in $M$ converging to an end $e$, and $B'_m,
m=1,2,\dots,$ is another such sequence converging to the same end, then if $n\geq 2$ there exists $g\in A_0$ such that $g(B_m)= B'_m$ for every $m$.
If in addition all the balls $B_m$ and $B'_\ell$ are pairwise disjoint, then there exists $g\in A_0$ such that $g(B_m)= B'_m$ and $g(B'_m)= B_m$ for every $m$. If
$n=1$, we can find $g\in A_0$ such that $g(B_m)= B'_m$ provided that the order of
the intervals $B_m$ and that of the intervals $B'_m$ on the line
is the same. \label{moveballs} \elm

\noindent{\bf Idea of the Proof.} 1. To obtain $g$, it suffices to
shrink $B$ and move it along $\gamma$ into the interior of $B'$
by a sequence of small isotopies.
Then, using the annulus theorem and a collar of $B'$, the image
of $B$ can be expanded to coincide with $B'$.

2. Using part (1) successively, the identity $I$ can be isotoped into $g$ so that $g(B_m)= B'_m$, and also
$g(B'_m)= B_m$ if the balls $B_m$ are disjoint from the balls $B'_\ell$, for all $m$ and $\ell$. Since the balls converge to the
end $e$, the successive isotopies can be done so that the points
moved are in smaller and smaller neighborhoods of $e$. Then
the composition of the isotopies gives an isotopy to the desired diffeomorphism $g$. \qed

\vskip 0.5cm

\noindent{\bf Proof of Proposition \ref{infballs}.}
Given a sequence of disjoint balls $B_m$ converging to an end $e$ of $M$, let $V_1\supset V_2\supset \dots$ be a base of product neighborhoods
of $e$. Let $F: M\times [0,1]\rightarrow M$, where
$F_t=F(\cdot,t)$ for $t\in [0,1]$, be an isotopy of the identity
map of $M$ to a diffeomorphism $F_1$ such that for all $m$,
$F_1(B_m)\subset \Int(V_1)$. This is easy
to carry out, using Lemma \ref{moveballs}, since each $V_k$ contains all
but finitely many of the original balls. Next let $F_t, t\in [1,2]$ be an
isotopy from $F_1$ to $F_2$ so that $F_2$ moves the balls $F_1(B_m)$ into
$V_2$ and such that the corresponding isotopy from $\I$ to $F_2
F_1^{-1}$ has support in $V_1$. When the isotopy from $F_{k-1}$ to
$F_k$ has been defined so that $F_k(B_m)\subset V_k$ for all $m$,
construct a further isotopy $F_t, \ t\in[k,k+1]$, from $F_k$ to $F_{k+1}$ so that $F_{k+1}$ moves the
balls $F_k(B_m)$ into $V_{k+1}$ and such that the corresponding isotopy
from $\I$ to $F_{k+1}F_k^{-1}$ has support in $V_k$. Thus we construct an
isotopy $F:M\times [0,\infty)\rightarrow M$ (straightening the
corners, if necessary), so that $F_k(B_m) \subset V_k$ for every
$m$ and when $x\in M\sm V_k,\ F(x,t)=F(x,k)$ for every $t\geq k$.

Now let $g$ be a diffeomorphism supported in $\cup_{m=1}^\infty
B_m$. We can define
an isotopy $G:M\times [0,1]\rightarrow M$ from $\I$ to $g$ by
setting $G(x,0)= x$ and $G(x,t) = F_{1/t} g F_{1/t}^{-1}(x)$ when $t\in (0,1]$. To see that $G$ is smooth along $M\times\{0\}$, note that $G(x,t) = F_{1/t} g F_{1/t}^{-1}(x)=x=G(x,0)$ when $x\in M\sm V_k$ and $1/t\geq k$. Thus $g\in A_0$. \qed

\vskip 0.5cm

The proof of Proposition \ref{infballs} can be refined to
show that every diffeomorphism of an open manifold $M$ with support in
a locally finite set of disjoint closed collared balls
(i.e., every compact set meets only finitely many of them)
is isotopic to the identity map of $M$.

\vskip 0.5cm

\noindent{\bf Proof of Proposition \ref{disjointballs}.}
Suppose that $h\in A_0=A_0(M)$ is a diffeomorphism whose support
meets every neighborhood of the end $e$ of $M$ and let
$\dim(M)\geq 2$. Choose a sequence of points $x_1,\dots,
x_m,\dots$ that converges to the end $e$ such that the points
$x_1,h(x_1),\dots,x_m,h(x_m),\dots,$ are distinct and form a
discrete set. Since the points $x_m$ converge to $e$ and $h$ is
isotopic to $I$, the points $h(x_m)$ must also converge to $e$.
Next
 choose closed collared balls $D_1,\dots,D_m,\dots$
with $x_m\in D_m$ so that the balls and their images $h(D_m)$ are all pairwise disjoint and they converge to the end $e$. Let $g\in A_0$
be a diffeomorphism with support inside a product neighborhood of the end $e$ that exchanges the balls $h(D_m)$ and $D_{m+1}$
for every $m$, by Lemma \ref{moveballs}. Then $f_1=(g^{-1}hg)h$ is a diffeomorphism in every normal subgroup $H$ that contains $h$ and
$f_1(D_m)=D_{m+2}$ for $m\geq 1$.
Now suppose that a sequence of disjoint balls $B_m$ converging
to the end $e$ is given, and let $g_1\in A_0$ be a diffeomorphism
such that $g_1(B_m)=D_{2m}$, again using Lemma \ref{moveballs}. Then $f=g_1^{-1}f_1g_1\in H$
satisfies $f(B_m)=g_1^{-1}f_1g_1(B_m)= g_1^{-1}f_1(D_{2m})=B_{m+1}$, as desired.

If $\dim(M)=1$, so that $M\approx {\mathbb R}$ with $e=+\infty$,
it is easy to show that $H$ contains a diffeomorphism $f_1$ that
is strictly increasing near $+\infty$. If $B_{m+1}$ lies to the
right of $B_m$ for every $m$, then $f_1$ can be conjugated in
$A_0$ to a diffeomorphism $f\in H$ that satisfies
$f(B_m)=B_{m+1}$.

It only remains to consider a
diffeomorphism $g$ of $M$ with $\supp(g)\subseteq
\cup_{m=1}^\infty B_m$ and show that $g\in H$. Using the same $f\in H$, define the diffeomorphism $\bar g$ to be
the identity outside $\cup_{m=2}^\infty B_m$ and to coincide with
$f\bar g g f^{-1}$ on $B_{m+1}$, once $\bar g$ has been defined on
$B_m$. (Note that $f^{-1}(B_{m+1})$ does not meet any ball
$B_k$ with $k\neq m$, so $\bar g$ is well-defined on $B_{m+1}$.)
This produces a well-defined diffeomorphism $\bar g$ of $M$, for the balls
$B_m$ are disjoint and locally finite so the local definitions
combine to give a global diffeomorphism. Furthermore $\bar g\in
A_0$ by Proposition \ref{infballs}. Now by construction $\bar g=f\bar g g f^{-1}$
on $\cup_{m=1}^\infty B_m$, and on the rest
of $M$ both coincide with the identity map. Thus $g=\bar
g^{-1}f^{-1}\bar g f$ belongs to $H$ since it is a
conjugate of $f^{-1}$ composed with $f$.
 \qed

\vskip 0.5cm

We remark that the proof of Proposition \ref{disjointballs} is
similar to steps in the proof of Lemma 2.7 of \cite{McD} and
related to results in \cite{L1}.

\vskip 0.5cm

\noindent{\bf Proof of Theorem \ref{thm2}.} Supposing the hypotheses, we claim that
$h$ must move points that converge to some end $e$ of $M$. To see
this, let $K_m$ be an increasing nested sequence of compact sets
such that $K_m\subset \Int(K_{m+1})$ for every $m$ and
$\cup_{m=1}^\infty K_m = M$. Choose a connected component $U_m$ of
$M\sm K_m$ for $m=1,2,\dots,$ so that the closure of $\supp(h)\cap
U_m$ is not compact, and so that $U_m\supset U_{m+1}$. Then the
sets $U_m$ form a neighborhood base of a unique end $e$, and $h$
moves points arbitrarily close to $e$.

Now if the given diffeomorphism $g$ has support in a ball, then
by Proposition \ref{disjointballs}, $g\in H$.
In general, if $g\in A_0(M)$ is any diffeomorphism isotopic to the identity
map by an isotopy with compact support, then by Corollary
\ref{factorization2} of the Factorization Lemma, $g$ is a
composition of diffeomorphisms supported in balls, each
of them belonging to $H$, so again $g\in H$. \qed

\vskip 0.5cm

\noindent{\bf Proof of Proposition \ref{nbdbase}.} Let $e$ be a
product end of a manifold $M$. Suppose that $h\in A_0\sm A_e$ so
that $\supp(h)$ meets every neighborhood of $e$, and let $H$ be a
normal subgroup of $A_0$ containing $h$. Take a product
neighborhood $V_1\cong N \times (0,1]$ of the end $e$ and identify $V_1$ with $N \times (0,1]$. Then the
sets $V_m= N \times (0,1/m]$, $m=1,2,\dots$, form a basis of
nested closed neighborhoods of $e$. Let $T:V_1\rightarrow V_2$ be
a diffeomorphism defined by setting $T(x,s)=(x,t(s))$ where
$t:(0,1]\rightarrow (0,1/2]$ satisfies $t(1/m)=1/(m+1)$ (for
example, $t(s)=\frac{s}{s+1}$). Take a finite set $B_1, B_2,
\dots, B_k$ of closed balls in $\Int(V_1)\sm V_4$ whose interiors
cover the compact set $V_2\sm \Int(V_3)$ and such that for
$j=1,\dots,k$ and $m>0$, $T^m(B_j)\cap B_j=\emptyset$. Now we may
conjugate the diffeomorphism $f\in H$ given by Proposition
\ref{disjointballs} by an element of $A_0$ to a diffeomorphism
$f_j\in H$ supported in $V_1$ so that $f_j^m(B_j)= T^m(B_j)$ for
$m=1,2,\dots$. Let $g_j$ be a diffeomorphism supported on
$\cup_{m=0}^\infty T^m(B_j)$ so that for each point $(x,s)\in
\Int(B_j)$, $g_j(x,s)=(x,s')$ for some $s'=s'(x,s)<s$ and such
that $g_j$ commutes with $T$. By Proposition \ref{disjointballs},
each $g_j\in H$.

The diffeomorphism $g=g_1\circ\dots\circ g_k$ preserves the first
coordinate of $N \times (0,1]$ and decreases the second coordinate
for every $(x,s)\in V_2\sm \Int(V_3)$. Since $V_2\sm \Int(V_3)$ is
compact, some power of $g$ satisfies $g^p(V_2)\subseteq V_3$. Set
$f_0=g^p$, an element of $H$ with support in $V_1$. Since $g^p$
commutes with $T$, it follows that $f_0(V_m)=f_0 T^{m-2}(V_2)=
T^{m-2}g^p(V_2)\subseteq T^{m-2}(V_3)=V_{m+1}$. The sets
$N_m:=f_0^m(V_2)\subseteq V_{m+2}$ for $m\geq 0$ form a neighborhood
base of $e$ consisting of nested product neighborhoods and
$f_0(N_m)=N_{m+1}$, so $f_0$ is the desired diffeomorphism. \qed

\vskip 0.5cm

\noindent{\bf Proof of Corollary \ref{prodsupp}.} Note that
any diffeomorphism with support in a product neighborhood of the
end $e$ is conjugate in $A_0$ to a diffeomorphism supported in
$N_1$, so it suffices to show that a diffeomorphism $g$
with $\supp(g)\subset N_1$ belongs to $H$. Define a diffeomorphism $\bar g$ with
support in $N_1$ recursively by setting
$\bar g = g^{-1}f_0\bar gf_0^{-1}$ on $f_0^m(N_0\sm \Int(N_1))=N_m\sm
\Int(N_{m+1})$ for $m=1,2,\dots$. This recursion formula also
holds on $M\sm N_1$, where both $g$ and $\bar g$ coincide with the
identity. Furthermore $\bar g$ is globally defined since
$\cap_{m=1}^\infty N_m =\emptyset$, and $\bar g\in A_0$ since its
support can be pushed off to infinity. Then the recursion formula
shows that $g=f_0\bar gf_0^{-1}\bar g^{-1}$, which is
an element of the normal subgroup $H$. \qed

\vskip 0.5cm

\noindent {\bf Proof of Theorem \ref{mainlemma}.} We use the
notation and hypotheses of the last two proofs. It suffices to
show that the normal subgroup $H$ of $A_0$ generated by $h$ and
$A_e$ is all of $A_0$, so consider $g\in A_0$ and an isotopy
$F:M\times [0,1]\rightarrow M$ from the identity map to $g$. Take
an integer $q$ sufficiently large so that $F(N_q\times
[0,1])\subset \Int(N_1)$. Let $C$ be the compact set $C=N_q
\sm\Int(N_{q+1})$ and take an integer $r$ large enough so that
$F(C\times [0,1])\cap N_r=\emptyset$. Let $\bar U_1,\dots, \bar
U_k$ be balls contained in $\Int(N_1)\sm N_r$ such that their
interiors $U_j$ cover $F(C\times [0,1])$. Applying the
Factorization Lemma \ref{factorization} to $C$ and the balls
$U_j$, we obtain finitely many diffeomorphisms $g_0=\I,g_1,\dots,
g_p\in A_0$, isotopies of $\I$ to $f_i=g_ig_{i-1}^{-1}$ for
$i=1,\dots,p$, each supported in one of the balls $U_j$, and an
isotopy of $\I$ to $gg_p^{-1}$ supported on $M\sm C$. By
Proposition \ref{infballs}, the diffeomorphism $g_p=f_p\dots f_1$,
which is supported on the compact union of the balls $\bar U_j$,
belongs to $H$.

Now $C=N_q \sm\Int(N_{q+1})$ separates $M\sm N_q$ from
$\Int(N_{q+1})$ and $gg_p^{-1}$ is isotopic to $\I$ by an isotopy
supported on $M\sm C$, so it can be factored $gg_p^{-1}=h_1h_2$,
where $h_1\in A_0$ is isotopic to $\I$ by an isotopy supported in
$M\sm N_q$ and $h_2\in A_0$ is isotopic to $\I$ by an isotopy
supported in $N_{q+1}$. Thus $h_1\in A_e\subset H$ by the
definition of $A_e$ and $h_2\in H$ by Corollary \ref{prodsupp}, so $g=h_1h_2g_p
\in H$. \qed

\end{document}